\documentclass[preprint,11pt,times]{elsarticle}
\usepackage{geometry}
\usepackage{amssymb,amsmath} 
\usepackage{graphicx,color,comment}
\usepackage{url}
\usepackage{mathrsfs}
\usepackage{subeqnarray,comment} 
\usepackage{graphicx}
\usepackage{placeins}
\usepackage{float}
\usepackage{array}
\usepackage{subfigure}
\usepackage{indentfirst}
\usepackage{multicol}
\usepackage{multirow}

\usepackage{amssymb}
\newtheorem{Theorem}{Theorem}[section]
\newtheorem{lemma}{Lemma}[section]
\newdefinition{rmk}{Remark}[section]
\newtheorem{corollary}{Corollary}[section]
\newproof{proof}{Proof}

\numberwithin{equation}{section}

\usepackage{amsmath}
\usepackage{lipsum}
\usepackage{amsfonts}
\usepackage{graphicx}
\usepackage{epstopdf}
\usepackage{algorithm}
\usepackage{algorithmic}
\ifpdf
  \DeclareGraphicsExtensions{.eps,.pdf,.png,.jpg}
\else
  \DeclareGraphicsExtensions{.eps}
\fi

\journal{Journal of Computational and Applied Mathematics}

\begin{document}

\begin{frontmatter}



\title{fOGA: An Orthogonal Greedy Algorithm for Fractional Laplacian Problems} 


\author[1]{Ruitong Shan%
}
\ead{shanrt22@mails.jlu.edu.cn}
\author[2]{Young Ju Lee}
\ead{yjlee@txstate.edu}
\author[1]{Jiwei Jia\corref{cor1}}
\ead{jiajiwei@jlu.edu.cn}
\cortext[cor1]{Corresponding author}

\affiliation[1]{organization={School of Mathematics},
addressline={Jilin University},
postcode={130012},
city={Changchun},
country={China}}
\affiliation[2]{organization={Department of Mathematics},
addressline={Texas State University},
city={San Marcos},
country={United States of
America}}
\begin{abstract}
We extend energy-based orthogonal greedy approximation to integral
fractional Poisson problems by combining OGA with discrete fractional
operators. In one dimension, order-dependent Gr\"unwald--Letnikov
stencils are paired with boundary-corrected ReLU atoms. On the unit disk,
a fractional centered-difference operator is constructed on a larger
periodic lattice and restricted to the interior, while boundary-weighted
ReLU$^2$ atoms impose the homogeneous exterior condition. The resulting
symmetric positive definite matrices provide a common energy framework
for normalized atom selection and orthogonal projection. We establish
the corresponding discrete projection properties and separate operator
discretization from approximation in the selected spaces. One- and
two-dimensional experiments assess the method through $L^2$, $H^1$,
and sampled maximum errors. The results demonstrate the applicability
of the OGA extension to nonlocal problems and show the influence of
boundary regularity on approximation, particularly in the $H^1$ error.

\end{abstract}




\end{frontmatter}



\section{Introduction}
Fractional-order equations have attracted considerable attention in recent decades because of their ability to describe a wide range of complex physical phenomena involving memory effects, hereditary properties, and nonlocal interactions. Fractional models and related nonlocal equations have found applications in quantum mechanics \cite{guo2015fractional,laskin2000fractional}, materials science \cite{bates2006some}, image processing \cite{gilboa2009nonlocal}, and anomalous diffusion and porous-medium flow \cite{lisini2018gradient,logvinova2004fractional}. Compared with classical integer-order differential operators, fractional operators provide a flexible framework for modeling anomalous diffusion, long-range interactions, and multiscale physical processes. Among them, the fractional Laplace operator plays a central role in the mathematical description of nonlocal diffusion and has become a fundamental operator in many fractional partial differential equations.

Despite its modeling advantages, the fractional Laplace operator poses substantial challenges for numerical computation. Owing to its intrinsic nonlocality and singular integral structure, the direct evaluation and discretization of this operator are generally much more involved than those of the classical Laplace operator. In practical computations, the value of the fractional Laplacian at a given point depends on the solution over the entire domain, and in some formulations even on its values outside the computational domain. This nonlocal feature often leads to dense stiffness matrices, increased computational cost, and additional difficulties in the treatment of boundary conditions. Therefore, the development of accurate, stable, and efficient numerical methods for fractional Laplace problems remains an important topic in scientific computing.

A variety of numerical methods have been proposed for equations involving the fractional Laplace operator, including finite difference and quadrature methods \cite{hao2021fractional,huang2014numerical}, finite element methods \cite{acosta2017fractional,ainsworth2018towards}, spectral methods \cite{hao2021sharp,hao2020optimal}, meshless pseudospectral methods \cite{burkardt2021unified,rosenfeld2019mesh}, and, more recently, deep-learning-based methods \cite{pang2019fpinns}. Each class of methods has its own advantages and limitations. Finite difference methods are relatively easy to implement and are well suited to structured grids, while finite element methods offer greater flexibility for complex geometries. Spectral methods can achieve high-order accuracy for sufficiently smooth solutions, but are often less convenient for irregular domains or nonsmooth solutions. Deep learning methods, by contrast, provide a promising mesh-free alternative and are particularly attractive for high-dimensional problems, although their theoretical foundation and numerical reliability still require further investigation.

In recent years, neural-network-based methods for partial differential equations have become an active research topic. One of their main advantages lies in their strong approximation capability for complex and high-dimensional function spaces. Among these approaches, physics-informed neural networks (PINNs) \cite{raissi2019physics} are particularly influential. By incorporating governing equations together with boundary or initial conditions into the loss function, PINNs provide a flexible framework for solving both forward and inverse problems. They have been applied successfully to fractional differential equations \cite{pang2019fpinns}, fluid dynamics \cite{bihlo2022physics,jin2021nsfnets}, and phase-field models \cite{mattey2022novel}, demonstrating considerable potential in scientific computing. In addition to PINNs, several other neural-network-based methods have also had a significant impact on the numerical solution of PDEs, including the Deep Ritz method \cite{yu2018deep}, weak adversarial networks \cite{zang2020weak}, the Deep Nitsche method \cite{liao2019deep}, and DeepONet \cite{lu2021learning}.

Nevertheless, despite their impressive empirical success, neural network methods remain much less understood from a rigorous numerical analysis perspective than classical schemes such as finite difference and finite element methods. The main difficulties arise from the nonlinear parameterization of neural networks, the nonconvexity of the associated optimization problems, and the essential role of training algorithms. Consequently, several fundamental questions, including convergence, stability, approximation properties, and generalization behavior, remain largely open.

A notable attempt to address these issues is the greedy approximation framework for shallow neural networks developed by Xu and collaborators, including the relaxed greedy algorithm (RGA) and the orthogonal greedy algorithm (OGA) \cite{siegel2021greedy}. Inspired by ideas from finite element methods and nonlinear approximation theory, this framework provides a more transparent description of the underlying approximation spaces and makes rigorous error estimates possible. Subsequent works have further developed this direction from both theoretical and computational viewpoints. In particular, greedy neural network methods have been extended to indefinite elliptic problems \cite{hong2025greedy}, and related studies on linearized shallow ReLU$^k$ networks have established integral representations of Sobolev spaces, optimal approximation rates, saturation results, and effective numerical methods for high-dimensional PDEs \cite{liu2025integral,mao2025sharp,mao2026solving}. Together, these results suggest that OGA-based shallow neural network methods offer a promising analytical framework for connecting neural network approximation with classical numerical analysis.

Motivated by these developments, we extend the energy-based OGA framework
to fractional Poisson equations through a discrete fractional operator.
The central construction connects three components: a positive definite
discretization of the nonlocal operator, a dictionary satisfying the
homogeneous exterior condition, and energy-normalized greedy selection
followed by projection. This makes the discrete fractional operator the
link between the governing equation and neural network approximation.

In one dimension, the directional representation leads to
order-dependent GL stencils and boundary-corrected ReLU atoms. On the
disk, a centered-difference Fourier symbol is raised to a fractional
power on a larger periodic lattice before restriction to the interior.
Boundary-weighted ReLU$^2$ atoms provide the two-dimensional trial space.
For both constructions, we establish positivity of the discrete matrices
and derive the associated projection and error-decomposition properties.
The resulting formulation adapts the OGA selection-and-projection
procedure to the nonlocal setting.

The numerical study uses consistent error measures across dimensions
and independent evaluation grids. The two disk benchmarks illustrate
the effect of boundary regularity: the first has a bounded gradient,
whereas the constant-forcing solution has a singular gradient.
Discrete-solve and quadrature-refinement checks complement the measured
errors and help interpret their dependence on the approximation space.

Section~\ref{sec:model} introduces the operators and their discretization.
Section~\ref{sec:oga} develops the OGA formulation and its discrete
properties. Section~\ref{sec:NUM} presents the numerical experiments,
followed by the conclusion.

\section{Governing equation and discretization}
\label{sec:model}
For the fractional problem, we use the integral fractional Laplacian
with homogeneous exterior data:
\begin{equation}
(-\Delta)^{\alpha/2}u=f\quad\hbox{in }\Omega,
\qquad u=0\quad\hbox{in }\mathbb R^d\setminus\Omega,
\qquad 0<\alpha<2.
\label{eq:exterior-problem}
\end{equation}
For sufficiently regular functions the operator is
\begin{equation}
(-\Delta)^{\alpha/2}u(\boldsymbol x)
=c_{d,\alpha}\operatorname{p.v.}\!\int_{\mathbb R^d}
\frac{u(\boldsymbol x)-u(\boldsymbol y)}{|\boldsymbol x-\boldsymbol y|^{d+\alpha}}
\,d\boldsymbol y,
\quad c_{d,\alpha}=\frac{2^\alpha\Gamma((d+\alpha)/2)}
{\pi^{d/2}|\Gamma(-\alpha/2)|}.
\end{equation}
The exterior values are part of the problem specification; boundary
values alone do not determine this nonlocal operator. This formulation
is distinct from the spectral fractional Laplacian on a bounded domain.
The one-dimensional analysis below treats $\Omega=(0,1)$ and
$\alpha\in(0,1)\cup(1,2)$, together with the classical case $\alpha=2$,
where $-u''=f$ and $u(0)=u(1)=0$. The one-dimensional GL implementation
does not include $\alpha=1$, for which that representation needs a
separate limiting formulation. The two-dimensional lattice construction
below has no such exclusion.

\subsection{Directional representation}
For $0<\alpha<2$ with $\alpha\ne1$, let $k=\lceil\alpha\rceil$. The directional derivative
and its relation to the fractional Laplacian are \cite{meerschaert2004vector}
\begin{equation}
D_{\boldsymbol\theta}^{\alpha}w(\boldsymbol x)
=\frac{(\boldsymbol\theta\cdot\nabla)^k}{\Gamma(k-\alpha)}
\int_0^\infty \xi^{k-\alpha-1}w(\boldsymbol x-\xi\boldsymbol\theta)\,d\xi,
\end{equation}
\begin{equation}
(-\Delta)^{\alpha/2}w(\boldsymbol x)
=C_{\alpha,d}\int_{S^{d-1}}D_{\boldsymbol\theta}^{\alpha}w(\boldsymbol x)
\,dS(\boldsymbol\theta),\qquad
C_{\alpha,d}=\frac{\Gamma((1-\alpha)/2)\Gamma((d+\alpha)/2)}
{2\pi^{(d+1)/2}}.
\end{equation}
Here $S^{d-1}$ is the unit sphere and $dS$ is surface measure. Gamma
values at negative noninteger arguments are understood by analytic
continuation. These representations are used where the derivatives and
integrals exist, or in their distributional interpretation.
For a zero extension from a convex domain, the ray integral can be
truncated at the backward distance
\begin{equation}
d(\boldsymbol x,\boldsymbol\theta,\Omega)
=\sup\{r>0:\boldsymbol x-\xi\boldsymbol\theta\in\Omega
\text{ for every }0\leq\xi<r\}.
\end{equation}
The dependence of this upper limit on $\boldsymbol x$ is retained under
differentiation. In a nonconvex domain a ray may leave and re-enter the
support, so the first boundary intersection alone need not suffice.

\subsection{One-dimensional GL discretization}
\label{sec:1d-discrete}
In one dimension the directions are $\theta=\pm1$, with backward
distances $x$ and $1-x$, respectively. The reflection identity gives
\begin{equation}
\Gamma\left(\frac{1+\alpha}{2}\right)
\Gamma\left(\frac{1-\alpha}{2}\right)
=\frac{\pi}{\cos(\pi\alpha/2)},\qquad
(-\Delta)^{\alpha/2}u
=\frac{D_+^\alpha u+D_-^\alpha u}{2\cos(\pi\alpha/2)}.
\end{equation}
The cosine is positive for $0<\alpha<1$ and negative for $1<\alpha<2$.
The GL coefficients and their partial sums are
\begin{equation}
B_k=(-1)^k\binom\alpha k,\qquad B_0=1,\qquad
B_k=\left(1-\frac{\alpha+1}{k}\right)B_{k-1},\quad k\geq1.
\label{eq:gl-coefficients}
\end{equation}
Let $m\geq2$, $h=1/m$, and $x_j=jh$. We distinguish the $m-1$ interior
grid values from the number $n$ of selected network atoms. With
$v_0=v_m=0$, define the shift
$s_\alpha=0$ for $0<\alpha<1$ and $s_\alpha=1$ for $1<\alpha<2$.
The discrete operator is
\begin{equation}
(A_h\boldsymbol v)_j=\gamma_{\alpha,h}
\left[\sum_{k=0}^{j+s_\alpha}B_kv_{j-k+s_\alpha}
+\sum_{k=0}^{m-j+s_\alpha}B_kv_{j+k-s_\alpha}\right],
\quad
\gamma_{\alpha,h}=\frac{1}{2\cos(\pi\alpha/2)h^\alpha}.
\label{eq:1d-gl}
\end{equation}
Here $j=1,\ldots,m-1$, and $\boldsymbol v=(v_1,\ldots,v_{m-1})^T$.
Both sums include $k=0$; their last terms are the prescribed zero
endpoint values. The choice of shift is essential for the positivity
properties established below.
The discrete problem is
\begin{equation}
A_h\boldsymbol u_h=\boldsymbol f_h,\qquad
\boldsymbol f_h=(f(x_1),\ldots,f(x_{m-1}))^T.
\label{eq:discrete-problem}
\end{equation}
For $1<\alpha<2$, the matrix entries are
\begin{equation}
(A_h)_{ij}=\gamma_{\alpha,h}
\begin{cases}
2B_1,&i=j,\\
B_0+B_2,&|i-j|=1,\\
B_{|i-j|+1},&|i-j|\geq2.
\end{cases}
\label{eq:shifted-matrix}
\end{equation}
For $\alpha=2$ we use $h^{-2}\operatorname{tridiag}(-1,2,-1)$.
The noninteger-order matrices are dense Toeplitz matrices, whereas the
classical matrix is tridiagonal.
The GL approximation is formally first order for sufficiently regular
zero extensions. This statement does not supply a uniform truncation
estimate at a nonsmooth boundary; no such regularity-independent estimate
is assumed in the analysis below.

\subsection{Positive definiteness of the one-dimensional matrices}
\label{sec:spd}
\begin{lemma}[Positive definiteness of the unshifted GL matrix]
\label{lem:unshifted-gl-spd}
Let $0<\alpha<1$, $m\geq2$, and $h=1/m$. After eliminating the zero
endpoint values, the unshifted GL operator in \eqref{eq:1d-gl} has the matrix
\begin{equation}
 (A_h)_{ij}=\gamma_{\alpha,h}
\begin{cases}
2,&i=j,\\
B_{|i-j|},&i\ne j,
\end{cases}
\qquad
\gamma_{\alpha,h}=\frac{1}{2\cos(\pi\alpha/2)h^\alpha}>0,
\label{eq:unshifted-gl-matrix}
\end{equation}
where $1\leq i,j\leq m-1$. Then $A_h$ is symmetric and strictly diagonally
dominant with positive diagonal entries, and hence is positive definite.
\end{lemma}

\begin{proof}
The GL coefficients satisfy
\[
B_0=1,\qquad B_1=-\alpha<0,\qquad
B_k=\left(1-\frac{\alpha+1}{k}\right)B_{k-1},\quad k\geq2.
\]
For $0<\alpha<1$, the multiplier is positive for every $k\geq2$,
so $B_k<0$ for all $k\geq1$. Moreover, their finite partial sums obey
\begin{equation}
S_j:=\sum_{k=0}^{j}B_k
=(-1)^j\binom{\alpha-1}{j}
=\prod_{k=1}^{j}\left(1-\frac{\alpha}{k}\right)>0,
\quad j\geq0,
\label{eq:gl-partial-sum}
\end{equation}
with $S_0=1$ and the empty product equal to one. The identity follows
from the generalized Pascal identity (or by induction).
For row $i$, the sum of the absolute values of the off-diagonal entries is
\[
\sum_{j\ne i}|(A_h)_{ij}|
=\gamma_{\alpha,h}
\left(-\sum_{k=1}^{i-1}B_k-\sum_{k=1}^{m-1-i}B_k\right).
\]
Consequently its diagonal-dominance margin is
\begin{equation}
d_i:=(A_h)_{ii}-\sum_{j\ne i}|(A_h)_{ij}|
=\gamma_{\alpha,h}\bigl(S_{i-1}+S_{m-1-i}\bigr)>0.
\label{eq:unshifted-gl-margin}
\end{equation}
The entries depend only on $|i-j|$, so the matrix is symmetric. Since
its off-diagonal entries are nonpositive, for any
$\boldsymbol v\in\mathbb R^{m-1}$ we have
\begin{equation}
\boldsymbol v^T A_h\boldsymbol v
=\sum_{i=1}^{m-1}d_i v_i^2
+\sum_{1\leq i<j\leq m-1}\bigl(-(A_h)_{ij}\bigr)(v_i-v_j)^2
\geq \left(\min_i d_i\right)\|\boldsymbol v\|_2^2.
\label{eq:unshifted-gl-positive-energy}
\end{equation}
The right-hand side is strictly positive when $\boldsymbol v\ne0$,
which proves the claim.
\end{proof}

\begin{lemma}[The shifted and classical cases]
\label{lem:shifted-spd}
For $1<\alpha<2$, the matrix \eqref{eq:shifted-matrix} is symmetric
positive definite. The same holds for the classical Dirichlet matrix.
\end{lemma}
\begin{proof}
For $1<\alpha<2$, $B_1=-\alpha<0$ and $B_k>0$ for $k\geq2$,
while $\gamma_{\alpha,h}<0$. Thus the diagonal is positive and every
off-diagonal entry is negative. The equalities in
\eqref{eq:gl-partial-sum} remain valid in this range, but now $S_j<0$
for $j\geq1$, while $S_0=1$. They yield $S_j\to0$ as $j\to\infty$.
Indeed, after its negative first factor,
the product has positive factors $1-\alpha/k$ for $k\geq2$ whose product
tends to zero. Hence $\sum_{k=0}^\infty B_k=0$, and
\[
2B_1+2(B_0+B_2)+2\sum_{k=3}^\infty B_k=0.
\]
The corresponding full-lattice row sums vanish. Each finite row omits
strictly negative off-diagonal entries on both tails, so its row sum is
strictly positive. It follows that $A_h$ is strictly diagonally dominant.
The quadratic-form argument in \eqref{eq:unshifted-gl-positive-energy}
therefore proves positive definiteness. For $\alpha=2$, with
$v_0=v_m=0$, the direct identity
\[
\boldsymbol v^TA_h\boldsymbol v
=h^{-2}\sum_{j=0}^{m-1}(v_{j+1}-v_j)^2>0
\quad\text{for }\boldsymbol v\ne0
\]
proves the result.
\end{proof}
If the shifted stencil were instead used for $0<\alpha<1$, its diagonal
would equal $2\gamma_{\alpha,h}B_1<0$. It could not define the positive
energy inner product used by our OGA. This is why the unshifted stencil
is used for the $\alpha=0.5$ experiments.

\subsection{Two-dimensional discretization on the disk}
\label{sec:2d-discrete}
Let $\Omega=\{\boldsymbol x\in\mathbb R^2:|\boldsymbol x|<1\}$.
We use a Cartesian lattice with $h=2/m$, where $m$ is an even integer,
and retain the strictly interior index set
\begin{equation}
\mathcal I_h=\{\boldsymbol i\in\mathbb Z^2:|h\boldsymbol i|<1\},
\qquad N_h=|\mathcal I_h|.
\end{equation}
The two-dimensional construction uses the fractional power of the
centered-difference Fourier symbol, in the setting of fractional
centered differences \cite{hao2021fractional}. It is specified below
including the finite periodic approximation used to compute its weights.

The centered-difference symbol on the full lattice is
\begin{equation}
\lambda_h(\boldsymbol\theta)
=\frac{4}{h^2}\left(\sin^2\frac{\theta_1}{2}
+\sin^2\frac{\theta_2}{2}\right).
\label{eq:2d-symbol}
\end{equation}
For $0<\alpha<2$, we use its fractional power $\lambda_h^{\alpha/2}$.
For fixed physical frequency $\boldsymbol\xi$,
$\lambda_h(h\boldsymbol\xi)^{\alpha/2}\to|\boldsymbol\xi|^\alpha$
as $h\to0$. This identifies the whole-space operator approximated by the lattice symbol.

Let $M>m$ be an even integer, and let $E$ inject a vector on
$\mathcal I_h$ into an $M\times M$ periodic lattice by setting all other
entries to zero. If $F_M$ is the unitary discrete Fourier transform and
$\boldsymbol k\in\{0,\ldots,M-1\}^2$, define
\begin{equation}
\begin{split}
B_{h,M}&=F_M^*\operatorname{diag}
 \left(\lambda_h(2\pi\boldsymbol k/M)^{\alpha/2}\right)F_M,\\
A_{h,M}&=E^TB_{h,M}E,\qquad
A_{h,M}\boldsymbol u_{h,M}=\boldsymbol f_h.
\end{split}
\label{eq:2d-restricted-matrix}
\end{equation}
The matrix is real because the symbol is real and even. Equivalently,
its entries are $w_{h,M}(\boldsymbol i-\boldsymbol j)$, where
\begin{equation}
w_{h,M}(\boldsymbol q)=\frac{1}{M^2}
\sum_{\boldsymbol k\in\{0,\ldots,M-1\}^2}
\lambda_h(2\pi\boldsymbol k/M)^{\alpha/2}
e^{2\pi\mathrm{i}\boldsymbol k\cdot\boldsymbol q/M}.
\end{equation}
The kernel is evaluated by an inverse FFT and restricted to the interior
indices. The fractional power is taken on the larger lattice
\emph{before} restriction. Taking the fractional power of a disk
Dirichlet matrix instead would approximate a different, spectral operator.
The construction in \eqref{eq:2d-restricted-matrix} is also defined at
$\alpha=1$ and reduces to the restricted five-point Laplacian at
$\alpha=2$, which is used as an implementation check.

The physical period is $L=Mh$. Finite $L$ introduces interactions with
periodic copies of the zero-extended disk data. Thus $A_{h,M}$ is a
finite-period approximation to a whole-space discretization, not an
exact representation of the infinite exterior. Grid refinement and
period enlargement address different errors. In particular, taking
$h\to0$ while keeping $L$ fixed does not remove periodic-image effects.
The experiments use $M=2^{\lceil\log_2(p_{\rm pad}m)\rceil}$ and check
the sensitivity to $p_{\rm pad}$ separately.

\subsection{Positive definiteness of the two-dimensional matrix}
\label{sec:2d-spd}
\begin{lemma}
\label{lem:2d-spd}
For $0<\alpha\leq2$ and $M>m$, the matrix $A_{h,M}$ in
\eqref{eq:2d-restricted-matrix} is symmetric positive definite.
\end{lemma}
\begin{proof}
All eigenvalues of $B_{h,M}$ are nonnegative. The symbol vanishes only
at the zero Fourier mode, so the nullspace of $B_{h,M}$ consists
precisely of constant lattice vectors. For a real interior vector
$\boldsymbol v$,
\begin{equation}
\boldsymbol v^TA_{h,M}\boldsymbol v
=\sum_{\boldsymbol k}
\lambda_h(2\pi\boldsymbol k/M)^{\alpha/2}
\bigl|(F_ME\boldsymbol v)_{\boldsymbol k}\bigr|^2\geq0.
\end{equation}
Equality implies that $E\boldsymbol v$ is constant on the entire
periodic lattice. Since $M>m$, the lattice contains nodes outside the
disk, where $E\boldsymbol v$ is zero. The constant is therefore zero,
and injectivity of $E$ gives $\boldsymbol v=0$. Symmetry follows from
the real even symbol.
\end{proof}
The positivity holds for each finite periodic lattice and provides the
energy inner product used below. Cholesky factorization verifies this
property in the implementation.

\section{Energy orthogonal greedy approximation}
\label{sec:oga}
We first describe the one-dimensional method using the matrices in
Section~\ref{sec:spd}. Section~\ref{sec:2d-oga} then gives the
two-dimensional dictionary, energy projection, and error decomposition
for the matrix established in Section~\ref{sec:2d-spd}.

\subsection{Boundary-corrected network dictionary}
For $t\in(0,1)$, define
\begin{equation}
g_t(x)=(x-t)_+-(1-t)x,\quad 0\leq x\leq1,
\qquad g_t(x)=0\quad\text{otherwise}.
\label{eq:boundary-atom}
\end{equation}
Each atom has zero endpoint values and belongs to $H_0^1(0,1)$;
its zero extension belongs to $H^1(\mathbb R)$ and hence to the
fractional energy space for $0<\alpha<2$.
Inside the interval its weak derivative is
$g_t'(x)=\boldsymbol1_{\{x>t\}}-(1-t)$ almost everywhere.
The derivative at the knot is immaterial to the Sobolev norm.
For the computations we use the finite dictionary
\begin{equation}
\mathcal D_{\mathcal T}=\{g_t:t\in\mathcal T\},\qquad
\mathcal T=\{\ell/1000:\ell=1,\ldots,999\}.
\end{equation}
The network after $n$ selections is
\begin{equation}
u_n(x)=\sum_{j=1}^n c_jg_{t_j}(x)
=\sum_{j=1}^n c_j(x-t_j)_+
-x\sum_{j=1}^n c_j(1-t_j),\quad x\in[0,1].
\end{equation}
Thus it is a shallow ReLU expansion with a linear boundary correction.
The reported count $n$ is the number of selected atoms, rather than a
parameter count matched to an arbitrary raw-ridge network.

\subsection{Discrete energy and residual selection}
Let $R_hv=(v(x_1),\ldots,v(x_{m-1}))^T$ denote nodal restriction.
On $\mathbb R^{m-1}$, define
\begin{equation}
a_h(\boldsymbol v,\boldsymbol w)=h\boldsymbol w^TA_h\boldsymbol v,
\qquad \|\boldsymbol v\|_{a_h}^2=a_h(\boldsymbol v,\boldsymbol v),
\qquad F_h(\boldsymbol w)=h\boldsymbol w^T\boldsymbol f_h.
\end{equation}
This is an inner product by Section~\ref{sec:spd}. It is not an inner
product on arbitrary continuous functions, since nodal restriction may
have a kernel. All projection statements below refer to nodal vectors.
For a nonzero candidate $\boldsymbol d_t=R_hg_t$, set
$\widehat{\boldsymbol d}_t=\boldsymbol d_t/\|\boldsymbol d_t\|_{a_h}$.
Starting with $\boldsymbol u_0=0$, choose
\begin{equation}
t_n\in\operatorname*{arg\,max}_{t\in\mathcal T_{n-1}}
\frac{|h\boldsymbol d_t^T(\boldsymbol f_h-A_h\boldsymbol u_{n-1})|}
{\|\boldsymbol d_t\|_{a_h}},
\label{eq:normalized-selection}
\end{equation}
where $\mathcal T_{n-1}$ excludes candidates with nodal vectors in the
current span. This is the Hilbert-space OGA correlation criterion
\cite{siegel2021greedy} applied to the normalized finite dictionary,
because $A_h\boldsymbol u_h=\boldsymbol f_h$ implies
\[
h\boldsymbol d_t^T(\boldsymbol f_h-A_h\boldsymbol u_{n-1})
=a_h(\boldsymbol u_h-\boldsymbol u_{n-1},\boldsymbol d_t).
\]
No knowledge of the exact PDE solution is needed for selection.

\subsection{Projection and finite-dimensional guarantees}
Let $G_n=[\boldsymbol d_{t_1},\ldots,\boldsymbol d_{t_n}]$ and
$V_n=\operatorname{range}(G_n)$. The coefficients solve
\begin{equation}
K_n\boldsymbol c_n=\boldsymbol b_n,
\qquad K_n=hG_n^TA_hG_n,\qquad
\boldsymbol b_n=hG_n^T\boldsymbol f_h,
\qquad \boldsymbol u_n=G_n\boldsymbol c_n.
\label{eq:oga-projection}
\end{equation}
This is an $n\times n$ system, independently of the grid size $m$.
\begin{Theorem}[Well-posed discrete projection]
\label{thm:projection}
If the columns of $G_n$ are linearly independent, $K_n$ is symmetric
positive definite and \eqref{eq:oga-projection} has a unique solution.
It satisfies
\begin{equation}
a_h(\boldsymbol u_h-\boldsymbol u_n,\boldsymbol v)=0
\quad(\boldsymbol v\in V_n),\qquad
\|\boldsymbol u_h-\boldsymbol u_n\|_{a_h}
=\min_{\boldsymbol v\in V_n}\|\boldsymbol u_h-\boldsymbol v\|_{a_h}.
\end{equation}
For nested selected spaces, the identity
\begin{equation}
\|\boldsymbol u_h-\boldsymbol u_{n-1}\|_{a_h}^2
=\|\boldsymbol u_h-\boldsymbol u_n\|_{a_h}^2
+\|\boldsymbol u_n-\boldsymbol u_{n-1}\|_{a_h}^2
\label{eq:pythagoras}
\end{equation}
holds. In exact arithmetic, if zero tolerance is used and selection
continues until no nonzero dictionary correlation remains, the resulting
vector is the projection of $\boldsymbol u_h$ onto the full dictionary
span $V_{\mathcal T}$.
\end{Theorem}
\begin{proof}
For $\boldsymbol c\ne0$, independence gives $G_n\boldsymbol c\ne0$,
so $\boldsymbol c^TK_n\boldsymbol c
=\|G_n\boldsymbol c\|_{a_h}^2>0$.
The coefficient equations give orthogonality. Expanding squared norms
and using that orthogonality proves the best-approximation property and
\eqref{eq:pythagoras}. If all correlations vanish, the residual is
orthogonal to every generator of $V_{\mathcal T}$, hence to its span.
Otherwise a nonzero correlation supplies an independent direction,
because the residual is already orthogonal to $V_n$. There can be at
most $\dim V_{\mathcal T}$ such additions.
\end{proof}
The final projection need not equal $\boldsymbol u_h$ unless the full
dictionary span contains it. In particular, the fixed 999-atom dictionary
need not span the entire interior nodal space when $m=4000$ or $8000$.

The implementation evaluates the finite candidate set, performs two
passes of energy Gram--Schmidt orthogonalization, and recovers the raw
atom coefficients by triangular back substitution. A candidate is
rejected when its orthogonalized energy norm is at most $10^{-10}$ times
its original energy norm. Selection stops at 64 atoms or when the largest
normalized correlation is below $10^{-13}$. These finite-precision rules
are implementation safeguards, not additional exact-arithmetic conclusions
of Theorem~\ref{thm:projection}.

\subsection{Error decomposition and scope of the analysis}
\label{sec:error-scope}
Assume that the exact solution admits point values and let
$\boldsymbol u^*=R_hu$. Its consistency residual is
$\boldsymbol\tau_h=A_h\boldsymbol u^*-\boldsymbol f_h$.
For a residual vector define the discrete dual norm by
\begin{equation}
\|\boldsymbol r\|_{a_h,*}^2
=h\boldsymbol r^TA_h^{-1}\boldsymbol r.
\end{equation}
Since $A_h(\boldsymbol u^*-\boldsymbol u_h)=\boldsymbol\tau_h$,
\begin{equation}
\|\boldsymbol u^*-\boldsymbol u_n\|_{a_h}
\leq\|\boldsymbol\tau_h\|_{a_h,*}
+\min_{\boldsymbol v\in V_n}\|\boldsymbol u_h-\boldsymbol v\|_{a_h}.
\label{eq:discrete-error-split}
\end{equation}
This separates the discretization error from approximation in the
selected nodal space. The residual after projection satisfies the exact
identity
\begin{equation}
\|\boldsymbol f_h-A_h\boldsymbol u_n\|_{a_h,*}
=\|\boldsymbol u_h-\boldsymbol u_n\|_{a_h}.
\end{equation}
The weighted Euclidean residual norm
$\sqrt h\|\boldsymbol f_h-A_h\boldsymbol u_n\|_2$ is a different
quantity and is not minimized by the energy projection.

For continuous error estimates, let $I_h$ be piecewise linear nodal
interpolation. For $X=L^2(0,1)$ or $H^1(0,1)$ whenever defined,
\begin{equation}
\begin{split}
\|u-u_n\|_X\leq{}&\|u-I_h\boldsymbol u^*\|_X
+\|I_h(\boldsymbol u^*-\boldsymbol u_h)\|_X\\
&+\|I_h(\boldsymbol u_h-\boldsymbol u_n)\|_X
+\|I_h\boldsymbol u_n-u_n\|_X.
\end{split}
\label{eq:continuous-split}
\end{equation}
In the reported experiments all dictionary knots are grid nodes, because
$m$ is a multiple of 1000. The last term then vanishes for the ReLU
expansion. The remaining terms connect the discrete projection with continuous
approximation through the regularity of the solution and the refinement
of the grid and dictionary.

\subsection{Two-dimensional dictionary and energy projection}
\label{sec:2d-oga}
On the disk, we impose the homogeneous exterior condition through the
trial functions
\begin{equation}
g_{\boldsymbol w,b}(\boldsymbol x)
=(1-|\boldsymbol x|^2)_+
\sigma(\boldsymbol w\cdot\boldsymbol x+b)^2,
\qquad \sigma(z)=\max(z,0),\quad |\boldsymbol w|=1.
\label{eq:2d-atoms}
\end{equation}
Each atom belongs to $H_0^1(\Omega)$, and its zero extension belongs to
$H^1(\mathbb R^2)$. Indeed, the interior gradient is bounded, the
function vanishes on the boundary, and the disk is a Lipschitz domain.
Hence the atoms are admissible in the fractional energy space for
$0<\alpha<2$. The continuous approximation is
$u_n=\sum_{j=1}^n c_jg_{\boldsymbol w_j,b_j}$.
It is a boundary-weighted shallow ReLU$^2$ expansion; $n$ denotes
selected atoms rather than the parameter count of an unconstrained network.

The finite dictionary uses
\begin{equation}
\begin{split}
\boldsymbol w_q&=(\cos(2\pi q/32),\sin(2\pi q/32)),
\quad q=0,\ldots,31,\\
b_\ell&=-0.8+1.8\ell/32,\quad \ell=0,\ldots,32.
\end{split}
\label{eq:2d-dictionary}
\end{equation}
There are 1056 candidate atoms. The lower bias cutoff avoids very narrow
caps near the boundary that are poorly resolved by the training lattice.
The same fixed dictionary is used for both benchmarks, without including
their exact solutions as candidate functions.

Let $R_h$ denote evaluation at the $N_h$ interior nodes. For a fixed
periodic lattice, write $A_h=A_{h,M}$ and let $\boldsymbol u_h$ solve
$A_h\boldsymbol u_h=\boldsymbol f_h$. We use the shorter notation
\begin{equation}
a_h(\boldsymbol v,\boldsymbol w)=h^2\boldsymbol w^TA_h\boldsymbol v,
\qquad \|\boldsymbol v\|_{a_h}^2=a_h(\boldsymbol v,\boldsymbol v),
\qquad F_h(\boldsymbol w)=h^2\boldsymbol w^T\boldsymbol f_h.
\end{equation}
The dependence on the period remains implicit in $A_h$ throughout the
two-dimensional projection and error analysis.

By Lemma~\ref{lem:2d-spd}, this is an inner product on nodal vectors.
It is not an inner product on all continuous trial functions, because
nodal restriction can have a nontrivial kernel. For a nonzero candidate
$\boldsymbol d=R_hg$, normalized residual selection maximizes
\begin{equation}
\frac{|h^2\boldsymbol d^T(\boldsymbol f_h-A_h\boldsymbol u_{n-1})|}
{\sqrt{h^2\boldsymbol d^TA_h\boldsymbol d}}
=\frac{|a_h(\boldsymbol u_h-\boldsymbol u_{n-1},\boldsymbol d)|}
{\|\boldsymbol d\|_{a_h}}.
\label{eq:2d-selection}
\end{equation}
Selection uses only the right-hand side and the discrete operator.
With $G_n$ formed from the selected nodal vectors, the coefficients satisfy
\begin{equation}
h^2G_n^TA_hG_n\boldsymbol c_n=h^2G_n^T\boldsymbol f_h,
\qquad \boldsymbol u_n=G_n\boldsymbol c_n=R_hu_n.
\label{eq:2d-projection}
\end{equation}

\begin{corollary}[Two-dimensional discrete projection]
\label{cor:2d-projection}
If the selected nodal vectors are linearly independent, the system
\eqref{eq:2d-projection} has a unique solution. For
$V_n=\operatorname{range}(G_n)$,
\begin{equation}
\|\boldsymbol u_h-\boldsymbol u_n\|_{a_h}
=\min_{\boldsymbol v\in V_n}
\|\boldsymbol u_h-\boldsymbol v\|_{a_h},
\end{equation}
and, for nested selected spaces,
\begin{equation}
\|\boldsymbol u_h-\boldsymbol u_{n-1}\|_{a_h}^2
=\|\boldsymbol u_h-\boldsymbol u_n\|_{a_h}^2
+\|\boldsymbol u_n-\boldsymbol u_{n-1}\|_{a_h}^2.
\end{equation}
In exact arithmetic, exhaustive nonzero-correlation selection terminates
at the energy projection onto the span of the finite nodal dictionary.
\end{corollary}
\begin{proof}
Lemma~\ref{lem:2d-spd} implies that $h^2G_n^TA_hG_n$ is positive definite
for independent columns. The coefficient equations give energy
orthogonality to $V_n$. The proof of Theorem~\ref{thm:projection}, with
weight $h^2$ in place of $h$, then gives the stated conclusions.
\end{proof}
The full dictionary span has dimension at most 1056, whereas the reported
training grid has 12849 nodes. The projection therefore need not recover
the full discrete solution. The monotonicity concerns the discrete energy distance; the continuous
errors are evaluated separately.

The implementation uses two passes of energy Gram--Schmidt and triangular
back substitution. Candidates with orthogonalized norm at most $10^{-9}$
times their original norm are rejected. Selection stops at 128 atoms or
when the largest normalized correlation is below $10^{-13}$.
The direct Cholesky solution of $A_h\boldsymbol u_h=\boldsymbol f_h$
serves as a discrete reference for projection checks.

\subsection{Two-dimensional error decomposition and regularity}
\label{sec:2d-error}
Assume that $u$ admits the required nodal values and write
$\boldsymbol u^*=R_hu$. Define the consistency residual and its associated discrete
dual norm by
\begin{equation}
\boldsymbol r_h=A_h\boldsymbol u^*-\boldsymbol f_h,
\qquad
\|\boldsymbol r\|_{a_h,*}^2
=h^2\boldsymbol r^TA_h^{-1}\boldsymbol r.
\end{equation}
\begin{Theorem}[Discrete error decomposition on the disk]
\label{thm:2d-error}
Under the independence assumption of Corollary~\ref{cor:2d-projection},
\begin{equation}
\|\boldsymbol u^*-\boldsymbol u_n\|_{a_h}
\leq\|\boldsymbol r_h\|_{a_h,*}
+\min_{\boldsymbol v\in V_n}\|\boldsymbol u_h-\boldsymbol v\|_{a_h}.
\label{eq:2d-error-split}
\end{equation}
Moreover,
\begin{equation}
\|\boldsymbol f_h-A_h\boldsymbol u_n\|_{a_h,*}
=\|\boldsymbol u_h-\boldsymbol u_n\|_{a_h}.
\end{equation}
\end{Theorem}
\begin{proof}
Since $A_h(\boldsymbol u^*-\boldsymbol u_h)
=\boldsymbol r_h$, the two corresponding primal and dual norms
are equal. The triangle inequality and the best-approximation property
give \eqref{eq:2d-error-split}. Applying the same norm identity to
$A_h(\boldsymbol u_h-\boldsymbol u_n)$ proves the residual statement.
\end{proof}

To distinguish lattice and finite-period effects, let $A_{h,\infty}$
be the interior restriction of the infinite-lattice convolution with kernel
\begin{equation}
w_{h,\infty}(\boldsymbol q)
=\frac{1}{(2\pi)^2}\int_{[-\pi,\pi]^2}
\lambda_h(\boldsymbol\theta)^{\alpha/2}
e^{\mathrm{i}\boldsymbol q\cdot\boldsymbol\theta}
\,d\boldsymbol\theta.
\end{equation}
The finite interior matrix is well defined by this integral. Algebraically,
\begin{equation}
\boldsymbol r_h
=\underbrace{(A_{h,\infty}\boldsymbol u^*-\boldsymbol f_h)}_
{\text{lattice consistency}}
+\underbrace{(A_h-A_{h,\infty})\boldsymbol u^*}_
{\text{finite-period error}}.
\label{eq:2d-consistency-split}
\end{equation}
Consequently, the first term of \eqref{eq:2d-error-split} is bounded
by the sum of the dual norms of these two residuals. Boundary regularity,
the lattice representation of the disk, and enlargement of the periodic
box must all be addressed to turn this decomposition into a convergence rate.

Continuous error control requires a reconstruction step as well.
For any specified linear reconstruction
$I_h:\mathbb R^{N_h}\to X$, where
$X=L^2(\Omega)$ or $H^1(\Omega)$ and $u,u_n\in X$,
\begin{equation}
\begin{split}
\|u-u_n\|_X\leq{}&\|u-I_h\boldsymbol u^*\|_X
+\|I_h(\boldsymbol u^*-\boldsymbol u_h)\|_X\\
&+\|I_h(\boldsymbol u_h-\boldsymbol u_n)\|_X
+\|I_h\boldsymbol u_n-u_n\|_X.
\end{split}
\label{eq:2d-continuous-split}
\end{equation}
This follows by adding and subtracting the three reconstructed vectors.
A reconstruction bound
$\|I_h\boldsymbol v\|_X\leq C_X(h,M)\|\boldsymbol v\|_{a_h}$
controls the middle terms, while the first and last terms describe
continuous reconstruction errors. In contrast to the aligned one-dimensional piecewise
linear setting, the last term need not vanish for the disk dictionary.
The continuous errors below are evaluated from the actual network,
not from a nodal interpolant.

Regularity is particularly relevant to the $H^1$ metric. For
$u=c(1-r^2)_+^p$, $c\ne0$, $p>0$, the change of variable $t=1-r^2$ gives
\begin{equation}
\int_\Omega|\nabla u|^2\,d\boldsymbol x
=4\pi c^2p^2\int_0^1(1-t)t^{2p-2}\,dt
=\frac{2\pi c^2p}{2p-1},\qquad p>\tfrac12.
\label{eq:disk-gradient-norm}
\end{equation}
The integral diverges for $0<p\leq1/2$. Thus the constant-forcing
solution with $p=\alpha/2$ belongs to $H^1(\Omega)$ exactly when
$\alpha>1$. Since every finite dictionary expansion belongs to $H^1$,
its error against that solution cannot have a finite $H^1$ norm when
$\alpha\leq1$. The two-dimensional tests therefore use $\alpha=1.2,1.6$.
Even in this range the gradient is singular near the boundary, so slow
$H^1$ approximation is compatible with decreasing discrete energy error.

\section{Numerical experiments}\label{sec:NUM}
\subsection{One-dimensional manufactured solution}
\label{sec:1d-results}
We solve \eqref{eq:exterior-problem} on $(0,1)$ with
$u(x)=x^3(1-x)^3$ inside the interval and zero outside.
For $\alpha\ne1$, the forcing in the interval is
\begin{equation}
f(x)=\frac{1}{2\cos(\pi\alpha/2)}
\sum_{p=3}^{6}a_p\frac{\Gamma(p+1)}{\Gamma(p+1-\alpha)}
\bigl[x^{p-\alpha}+(1-x)^{p-\alpha}\bigr],
\quad (a_3,a_4,a_5,a_6)=(1,-3,3,-1).
\end{equation}
For $\alpha=2$ this reduces to
$f(x)=-6x+36x^2-60x^3+30x^4$.
The exact solution is smooth on the interval and belongs to $H_0^1(0,1)$,
so all three reported error measures are well defined.
The finite dictionary, normalized selection, and projection are precisely
those of Section~\ref{sec:oga}. The fractional runs use $m=4000$ and
$8000$, and the classical comparison uses $m=4000$.

The $L^2$ and $H^1$ errors are approximated by
\begin{equation}
E_{2,h}=\left(\sum_q\eta_q|e(z_q)|^2\right)^{1/2},\qquad
E_{1,h}=\left(\sum_q\eta_q\bigl[|e(z_q)|^2+|e'(z_q)|^2\bigr]\right)^{1/2},
\quad e=u-u_n.
\end{equation}
The nodes and physical weights $(z_q,\eta_q)$ are obtained from four-point
Gauss rules on the partition containing $j/1000$ and all dictionary knots.
An independent eight-point rule with breakpoints at all selected knots
agreed to absolute tolerances $10^{-10}$ in $L^2$ and $10^{-9}$ in $H^1$
for all tabulated rows. The table headings use continuous norm notation
for these verified quadrature approximations.
The sampled maximum error is
\begin{equation}
\|u-u_n\|_{\infty,\mathrm{test}}=\max_{x\in X_{\mathrm{test}}}|u(x)-u_n(x)|,
\end{equation}
where $X_{\mathrm{test}}$ contains 20001 uniformly spaced points and all
dictionary knots. This is not a certified continuous supremum.
For each error measure, the tabulated empirical order is
$p_n=\log_2(E_{n/2}/E_n)$ at $n=4,8,16,32,64$; the first row is marked
by a dash. These orders refer to increasing atom counts on a fixed grid.

\subsubsection{Classical case: $\alpha=2$}
The classical run uses the tridiagonal Dirichlet matrix and the same
boundary-corrected atoms as the fractional runs. The decrease of the
solution errors should not
be interpreted as convergence of a strong-form PINN loss: a finite ReLU
expansion is piecewise linear and its second derivative has singular
contributions at the knots.
\begin{table}[H]
\centering
\caption{Classical case $\alpha=2$, $m=4000$ subintervals (3999 interior training points), using the same boundary-corrected dictionary and error evaluation as the fractional runs.}
\label{tab:classical-new}
\small
\renewcommand\arraystretch{1.3}
\begin{tabular}{r rr rr rr}
\hline
$n$ & $\|u-u_n\|_{L^2}$ & order & $\|u-u_n\|_{H^1}$ & order & $\|u-u_n\|_{\infty,\mathrm{test}}$ & order \\ \hline
2 & 1.32608e-03 & -- & 1.58574e-02 & -- & 2.62967e-03 & -- \\
4 & 8.09797e-04 & 0.712 & 1.14304e-02 & 0.472 & 1.95945e-03 & 0.424 \\
8 & 2.33932e-04 & 1.791 & 6.89076e-03 & 0.730 & 5.55251e-04 & 1.819 \\
16 & 8.48165e-05 & 1.464 & 3.77179e-03 & 0.869 & 1.78935e-04 & 1.634 \\
32 & 2.05521e-05 & 2.045 & 1.93597e-03 & 0.962 & 4.60546e-05 & 1.958 \\
64 & 5.14270e-06 & 1.999 & 9.81154e-04 & 0.981 & 1.14480e-05 & 2.008 \\
\hline
\end{tabular}
\end{table}

\subsubsection{Fractional case: $\alpha=1.5$}
We use the shifted GL matrix \eqref{eq:shifted-matrix}. Its positive
definiteness and the well-posedness of the coefficient projection follow
from Lemma~\ref{lem:shifted-spd} and Theorem~\ref{thm:projection}.
\begin{table}[H]
\centering
\caption{One-dimensional fractional problem with $\alpha=1.5$ and $m=4000$ subintervals (3999 interior training points). Orders are measured as the number $n$ of selected atoms doubles; $\|u-u_n\|_{\infty,\mathrm{test}}$ is the sampled maximum error.}
\label{tab:alpha1p5_m4000}
\renewcommand\arraystretch{1.3}
\small
\begin{tabular}{r rr rr rr}
\hline
$n$ & $\|u-u_n\|_{L^2}$ & order & $\|u-u_n\|_{H^1}$ & order & $\|u-u_n\|_{\infty,\mathrm{test}}$ & order \\ \hline
2 & 1.27817e-03 & -- & 1.60554e-02 & -- & 2.93243e-03 & -- \\
4 & 5.92780e-04 & 1.109 & 1.20813e-02 & 0.410 & 1.38919e-03 & 1.078 \\
8 & 1.49695e-04 & 1.985 & 7.56866e-03 & 0.675 & 3.70790e-04 & 1.906 \\
16 & 4.92223e-05 & 1.605 & 4.01383e-03 & 0.915 & 1.06995e-04 & 1.793 \\
32 & 1.23560e-05 & 1.994 & 2.09517e-03 & 0.938 & 3.78746e-05 & 1.498 \\
64 & 3.95328e-06 & 1.644 & 1.05336e-03 & 0.992 & 1.12002e-05 & 1.758 \\
\hline
\end{tabular}
\end{table}
\begin{table}[H]
\centering
\caption{One-dimensional fractional problem with $\alpha=1.5$ and $m=8000$ subintervals (7999 interior training points). Orders are measured as the number $n$ of selected atoms doubles; $\|u-u_n\|_{\infty,\mathrm{test}}$ is the sampled maximum error.}
\label{tab:alpha1p5_m8000}
\renewcommand\arraystretch{1.3}
\small
\begin{tabular}{r rr rr rr}
\hline
$n$ & $\|u-u_n\|_{L^2}$ & order & $\|u-u_n\|_{H^1}$ & order & $\|u-u_n\|_{\infty,\mathrm{test}}$ & order \\ \hline
2 & 1.27845e-03 & -- & 1.60566e-02 & -- & 2.93344e-03 & -- \\
4 & 5.92188e-04 & 1.110 & 1.20838e-02 & 0.410 & 1.38717e-03 & 1.080 \\
8 & 1.49140e-04 & 1.989 & 7.57017e-03 & 0.675 & 3.73620e-04 & 1.893 \\
16 & 4.88434e-05 & 1.610 & 4.01468e-03 & 0.915 & 1.07228e-04 & 1.801 \\
32 & 1.18919e-05 & 2.038 & 2.09561e-03 & 0.938 & 3.63244e-05 & 1.562 \\
64 & 3.14918e-06 & 1.917 & 1.05354e-03 & 0.992 & 9.55309e-06 & 1.927 \\
\hline
\end{tabular}
\end{table}
For $m=8000$, the final doubling from 32 to 64 atoms gives empirical
orders of 1.917 in $L^2$, 0.992 in $H^1$, and 1.927 in the sampled
maximum error. For $m=4000$, the corresponding $L^2$ order is 1.644.
The difference reflects the combined influence of discretization and
selected-space approximation on the measured error.

\subsubsection{Fractional case: $\alpha=0.5$}
We use the unshifted GL matrix \eqref{eq:unshifted-gl-matrix}, whose
positive definiteness is established in Lemma~\ref{lem:unshifted-gl-spd}.
All dictionary and error-evaluation settings are the same as above.
\begin{table}[H]
\centering
\caption{One-dimensional fractional problem with $\alpha=0.5$ and $m=4000$ subintervals (3999 interior training points). Orders are measured as the number $n$ of selected atoms doubles; $\|u-u_n\|_{\infty,\mathrm{test}}$ is the sampled maximum error.}
\label{tab:alpha0p5_m4000}
\renewcommand\arraystretch{1.3}
\small
\begin{tabular}{r rr rr rr}
\hline
$n$ & $\|u-u_n\|_{L^2}$ & order & $\|u-u_n\|_{H^1}$ & order & $\|u-u_n\|_{\infty,\mathrm{test}}$ & order \\ \hline
2 & 1.27896e-03 & -- & 1.67081e-02 & -- & 2.90772e-03 & -- \\
4 & 5.64376e-04 & 1.180 & 1.49551e-02 & 0.160 & 1.71038e-03 & 0.766 \\
8 & 1.33905e-04 & 2.075 & 7.43588e-03 & 1.008 & 4.10654e-04 & 2.058 \\
16 & 4.44291e-05 & 1.592 & 4.58819e-03 & 0.697 & 1.77395e-04 & 1.211 \\
32 & 1.12367e-05 & 1.983 & 2.26281e-03 & 1.020 & 4.87124e-05 & 1.865 \\
64 & 3.49030e-06 & 1.687 & 1.12855e-03 & 1.004 & 1.12971e-05 & 2.108 \\
\hline
\end{tabular}
\end{table}
\begin{table}[H]
\centering
\caption{One-dimensional fractional problem with $\alpha=0.5$ and $m=8000$ subintervals (7999 interior training points). Orders are measured as the number $n$ of selected atoms doubles; $\|u-u_n\|_{\infty,\mathrm{test}}$ is the sampled maximum error.}
\label{tab:alpha0p5_m8000}
\renewcommand\arraystretch{1.3}
\small
\begin{tabular}{r rr rr rr}
\hline
$n$ & $\|u-u_n\|_{L^2}$ & order & $\|u-u_n\|_{H^1}$ & order & $\|u-u_n\|_{\infty,\mathrm{test}}$ & order \\ \hline
2 & 1.27911e-03 & -- & 1.67096e-02 & -- & 2.90833e-03 & -- \\
4 & 5.64271e-04 & 1.181 & 1.49589e-02 & 0.160 & 1.71288e-03 & 0.764 \\
8 & 1.33843e-04 & 2.076 & 7.43802e-03 & 1.008 & 4.12893e-04 & 2.053 \\
16 & 4.43846e-05 & 1.592 & 4.59369e-03 & 0.695 & 1.79546e-04 & 1.201 \\
32 & 1.10580e-05 & 2.005 & 2.26418e-03 & 1.021 & 5.07518e-05 & 1.823 \\
64 & 2.89154e-06 & 1.935 & 1.12906e-03 & 1.004 & 1.14917e-05 & 2.143 \\
\hline
\end{tabular}
\end{table}
For $m=8000$, the final empirical orders are 1.935 in $L^2$, 1.004 in
$H^1$, and 2.143 in the sampled maximum error. These results describe
the smooth manufactured solution with a fixed dictionary and up to
64 selected atoms. The disk experiments below examine the effect of
fractional boundary singularities.

\FloatBarrier
\subsection{Two-dimensional disk benchmarks}
\label{sec:2d-results}
We consider two radial solutions of \eqref{eq:exterior-problem} on the
unit disk. The exact solution--forcing pairs below follow from the
fractional Laplacian formulas for $(1-|\boldsymbol x|^2)_+^p$
in \cite{dyda2012fractional}, with $p=1+\alpha/2$ and $p=\alpha/2$,
respectively, using the normalization in Section~\ref{sec:model}.
Both tests use $\alpha=1.2$ and $1.6$, the lattice operator
of Section~\ref{sec:2d-discrete}, and the dictionary in
\eqref{eq:2d-dictionary}. The training grid has $m=128$, $h=1/64$ and
$N_h=12849$ strictly interior nodes. With $p_{\rm pad}=16$, the FFT
side length is $M=2048$ and the physical period is $L=32$.
The exterior condition is built into the atoms, so no boundary penalty
or separate boundary training points are used. All four runs reach
128 selected atoms.

The error measures follow the one-dimensional
conventions in Section~\ref{sec:1d-results}. We report
$\|u-u_n\|_{L^2(\Omega)}$, $\|u-u_n\|_{H^1(\Omega)}$, and the sampled
maximum error $\|u-u_n\|_{\infty,\mathrm{test}}$, using independent quadrature
and test points. The $H^1$ error includes both the function and gradient
terms, with boundary-adapted quadrature for the singular solution in
Example~2. The tables and convergence curves start at $n=16$ and use
$n=16,32,64,128$, with the two fractional orders presented side by side
for each example.

\subsubsection{Example 1: a solution with bounded gradient}
Set $s=\alpha/2$. The exact solution and interior forcing are
\begin{equation}
\begin{split}
u(\boldsymbol x)&=(1-|\boldsymbol x|^2)_+^{1+s},\\
f(\boldsymbol x)&=2^\alpha\Gamma(2+s)\Gamma(1+s)
\bigl[1-(1+s)|\boldsymbol x|^2\bigr],\quad \boldsymbol x\in\Omega.
\end{split}
\label{eq:2d-example1}
\end{equation}
The positive part enforces zero exterior data. The gradient is bounded
and vanishes at the boundary.
\begin{table}[H]
\centering
\caption{Errors for Example 1 with $m=128$ (12849 interior training points).}
\label{tab:2d-ex1-errors}
\small
\renewcommand\arraystretch{1.3}
\setlength{\tabcolsep}{4pt}
\begin{tabular}{r rrr rrr}
\hline
& \multicolumn{3}{c}{$\alpha=1.2$} & \multicolumn{3}{c}{$\alpha=1.6$} \\
\cline{2-4}\cline{5-7}
$n$ & $\|u-u_n\|_{L^2}$ & $\|u-u_n\|_{H^1}$ & $\|u-u_n\|_{\infty,\mathrm{test}}$
& $\|u-u_n\|_{L^2}$ & $\|u-u_n\|_{H^1}$ & $\|u-u_n\|_{\infty,\mathrm{test}}$ \\ \hline
16 & 1.21460e-02 & 1.32633e-01 & 2.38626e-02 & 1.29050e-02 & 1.19497e-01 & 3.02369e-02 \\
32 & 3.20708e-03 & 5.27412e-02 & 7.06971e-03 & 1.76502e-03 & 3.33991e-02 & 4.41231e-03 \\
64 & 1.33054e-03 & 3.44479e-02 & 2.95268e-03 & 3.96972e-04 & 1.38359e-02 & 9.48245e-04 \\
128 & 4.88905e-04 & 2.62785e-02 & 8.61136e-04 & 2.89230e-04 & 1.12471e-02 & 6.89547e-04 \\
\hline
\end{tabular}
\end{table}

\begin{figure}[htbp]
\centering
\includegraphics[width=0.40\textwidth]{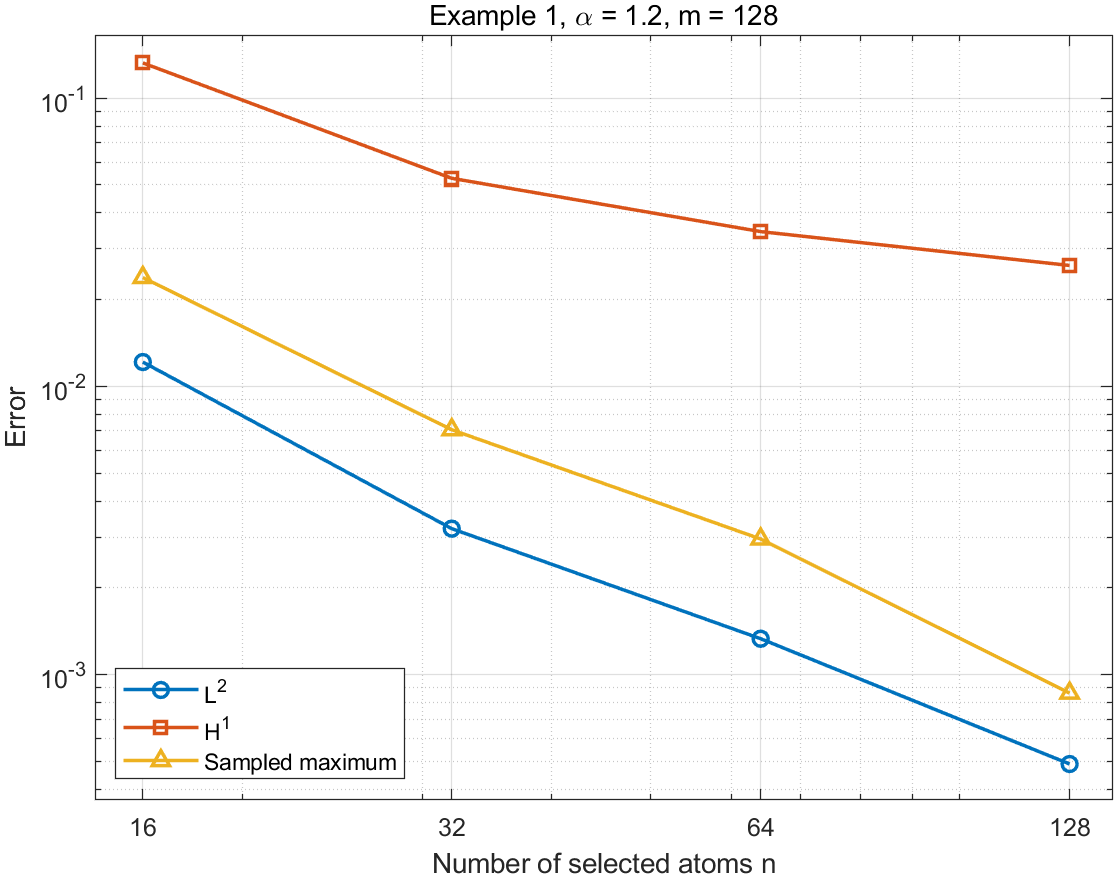}
\includegraphics[width=0.40\textwidth]{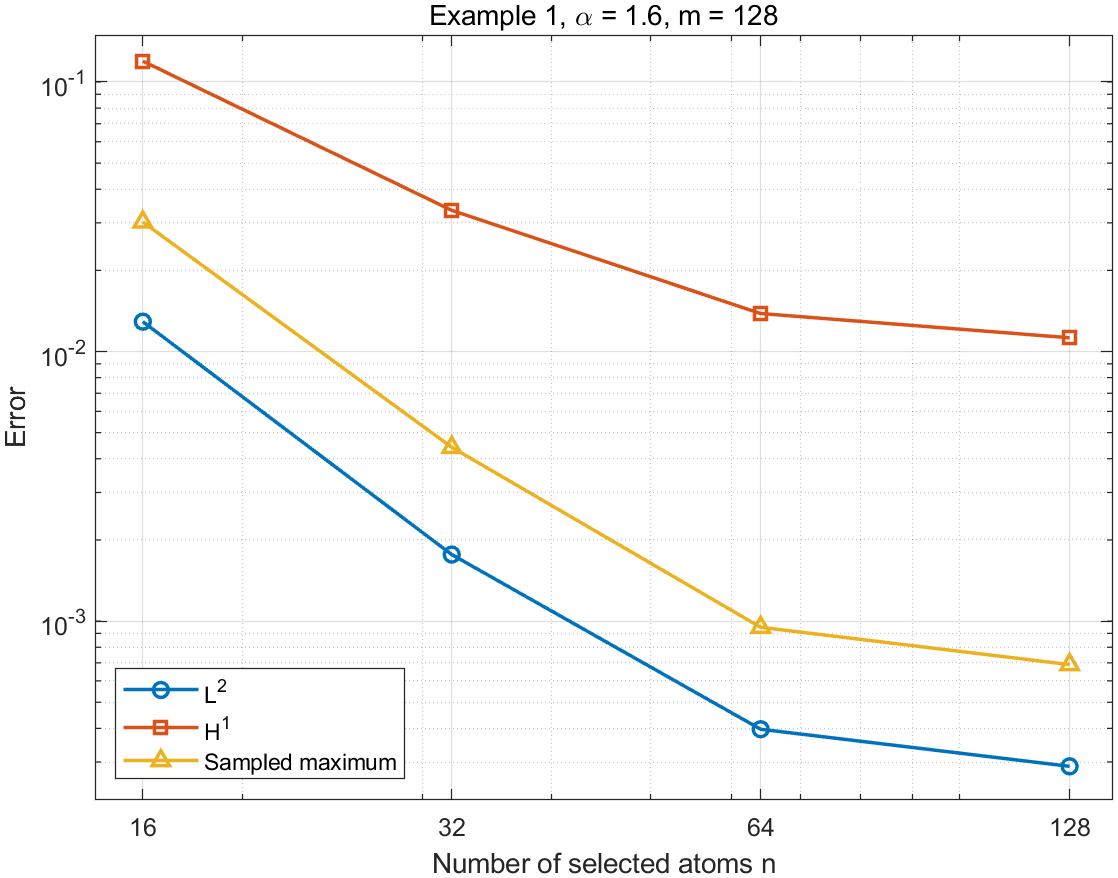}
\caption{Example 1: continuous error estimates against selected-atom count,
with $\alpha=1.2$ (left) and $\alpha=1.6$ (right). The maximum error is sampled.}
\label{fig:2d-ex1-convergence}
\end{figure}

\begin{figure}[htbp]
\centering
\includegraphics[width=0.84\textwidth]{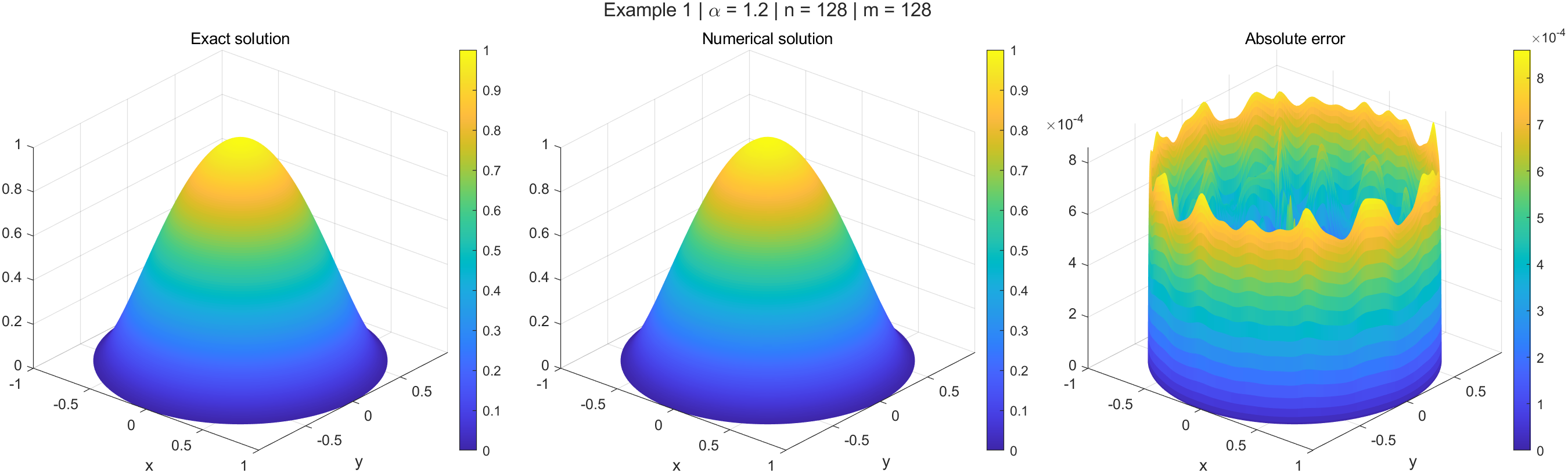}
\caption{Example 1, $\alpha=1.2$, $n=128$: exact solution (left),
network approximation (middle), and absolute error (right).
The first two panels share vertical and color scales.}
\label{fig:2d-ex1-a1p2-solutions}
\end{figure}

\begin{figure}[htbp]
\centering
\includegraphics[width=0.84\textwidth]{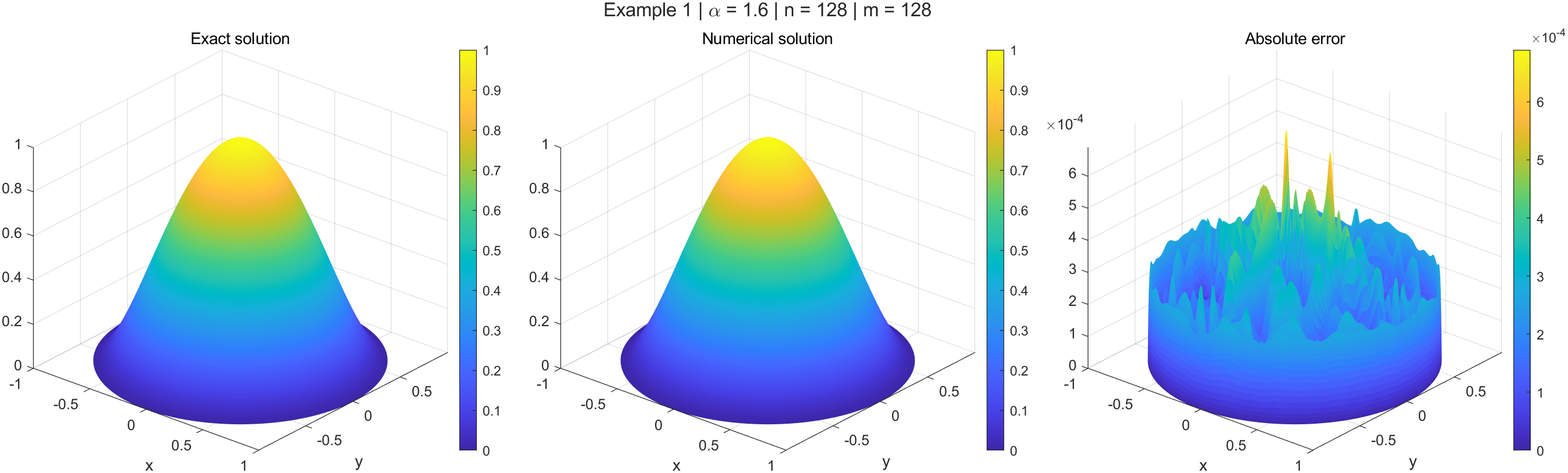}
\caption{Example 1, $\alpha=1.6$, $n=128$: exact solution (left),
network approximation (middle), and absolute error (right).
The first two panels share vertical and color scales.}
\label{fig:2d-ex1-a1p6-solutions}
\end{figure}

Table~\ref{tab:2d-ex1-errors} shows decreasing
errors at all reported atom counts. At $n=128$, the $L^2$ errors are
$4.889049\times10^{-4}$ and $2.892299\times10^{-4}$ for
$\alpha=1.2$ and $1.6$, respectively. The corresponding $H^1$ errors
are $2.627850\times10^{-2}$ and $1.124712\times10^{-2}$.
The gradient error improves more slowly at larger atom counts,
as seen in the last two rows. Figure~\ref{fig:2d-ex1-convergence} displays this behavior.
The solution comparisons in Figures~\ref{fig:2d-ex1-a1p2-solutions}
and \ref{fig:2d-ex1-a1p6-solutions} resolve the radial profile while also
showing the smaller angular variations of the approximation error.

\FloatBarrier
\subsubsection{Example 2: constant forcing and a boundary singularity}
For constant interior forcing, the solution is
\begin{equation}
f(\boldsymbol x)=1\quad(\boldsymbol x\in\Omega),\qquad
u(\boldsymbol x)=\frac{(1-|\boldsymbol x|^2)_+^s}
{2^\alpha\Gamma(1+s)^2},\quad s=\alpha/2.
\label{eq:2d-example2}
\end{equation}
Unlike Example~1, this solution has an unbounded boundary gradient for
both tested orders. Its $H^1$ norm is nevertheless finite because
$\alpha>1$, as established by \eqref{eq:disk-gradient-norm}.
\begin{table}[H]
\centering
\caption{Errors for Example 2 with $m=128$ (12849 interior training points).}
\label{tab:2d-ex2-errors}
\small
\renewcommand\arraystretch{1.3}
\setlength{\tabcolsep}{4pt}
\begin{tabular}{r rrr rrr}
\hline
& \multicolumn{3}{c}{$\alpha=1.2$} & \multicolumn{3}{c}{$\alpha=1.6$} \\
\cline{2-4}\cline{5-7}
$n$ & $\|u-u_n\|_{L^2}$ & $\|u-u_n\|_{H^1}$ & $\|u-u_n\|_{\infty,\mathrm{test}}$
& $\|u-u_n\|_{L^2}$ & $\|u-u_n\|_{H^1}$ & $\|u-u_n\|_{\infty,\mathrm{test}}$ \\ \hline
16 & 2.40216e-02 & 1.28325e+00 & 3.53654e-02 & 6.65945e-03 & 1.49355e-01 & 1.20851e-02 \\
32 & 9.38926e-03 & 1.26044e+00 & 1.81485e-02 & 2.54236e-03 & 1.32866e-01 & 3.93034e-03 \\
64 & 7.60828e-03 & 1.25357e+00 & 1.75644e-02 & 1.67514e-03 & 1.30318e-01 & 3.30387e-03 \\
128 & 7.18009e-03 & 1.25024e+00 & 1.67890e-02 & 1.58864e-03 & 1.29420e-01 & 2.75509e-03 \\
\hline
\end{tabular}
\end{table}

Table~\ref{tab:2d-ex2-errors} shows a smaller reduction in the $H^1$
error than in the other two error measures. For $\alpha=1.2$, the
$L^2$ error decreases from $2.402162\times10^{-2}$ at $n=16$ to
$7.180091\times10^{-3}$ at $n=128$, whereas the final $H^1$ error is
about 1.25, with little change between 64 and 128 atoms.
For $\alpha=1.6$, the final $L^2$ and $H^1$ errors are
$1.588641\times10^{-3}$ and $1.294197\times10^{-1}$, respectively.
The last two rows again show that the $H^1$ error is approaching a plateau
at the present grid and dictionary resolutions. 

\begin{figure}[htbp]
\centering
\includegraphics[width=0.40\textwidth]{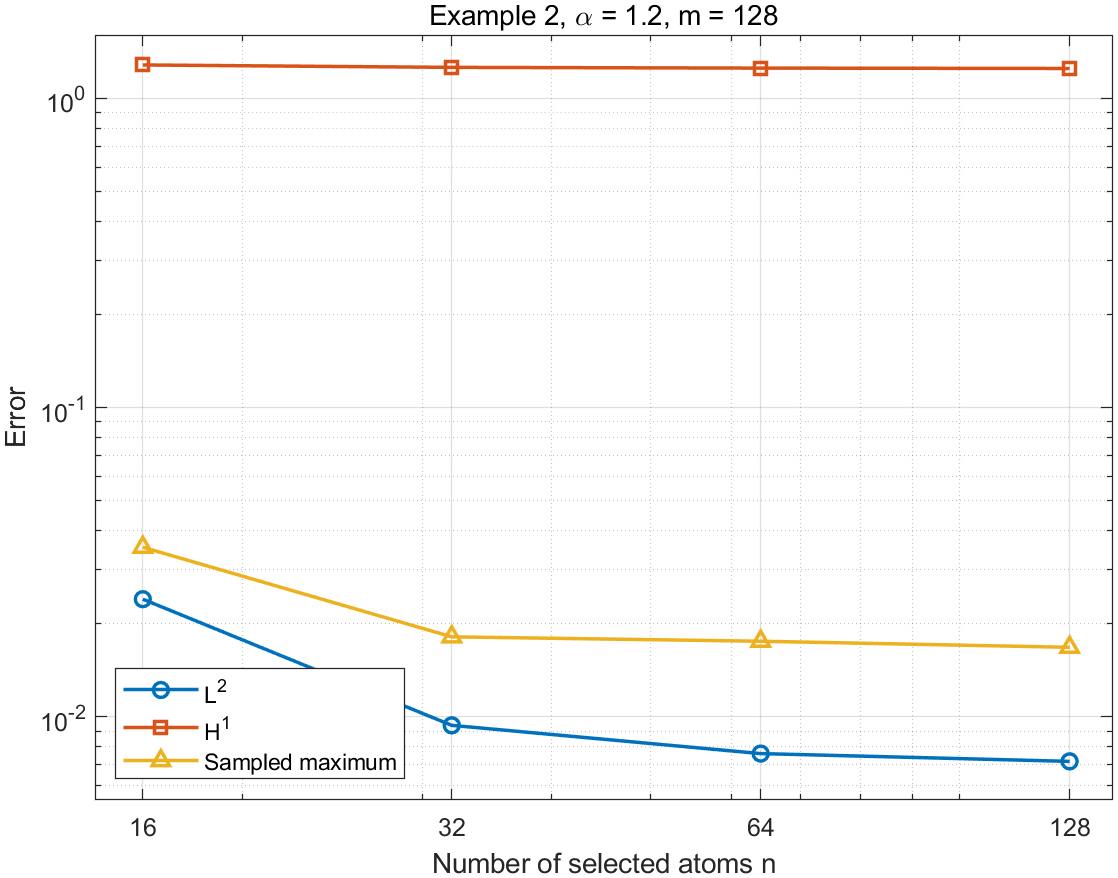}
\includegraphics[width=0.40\textwidth]{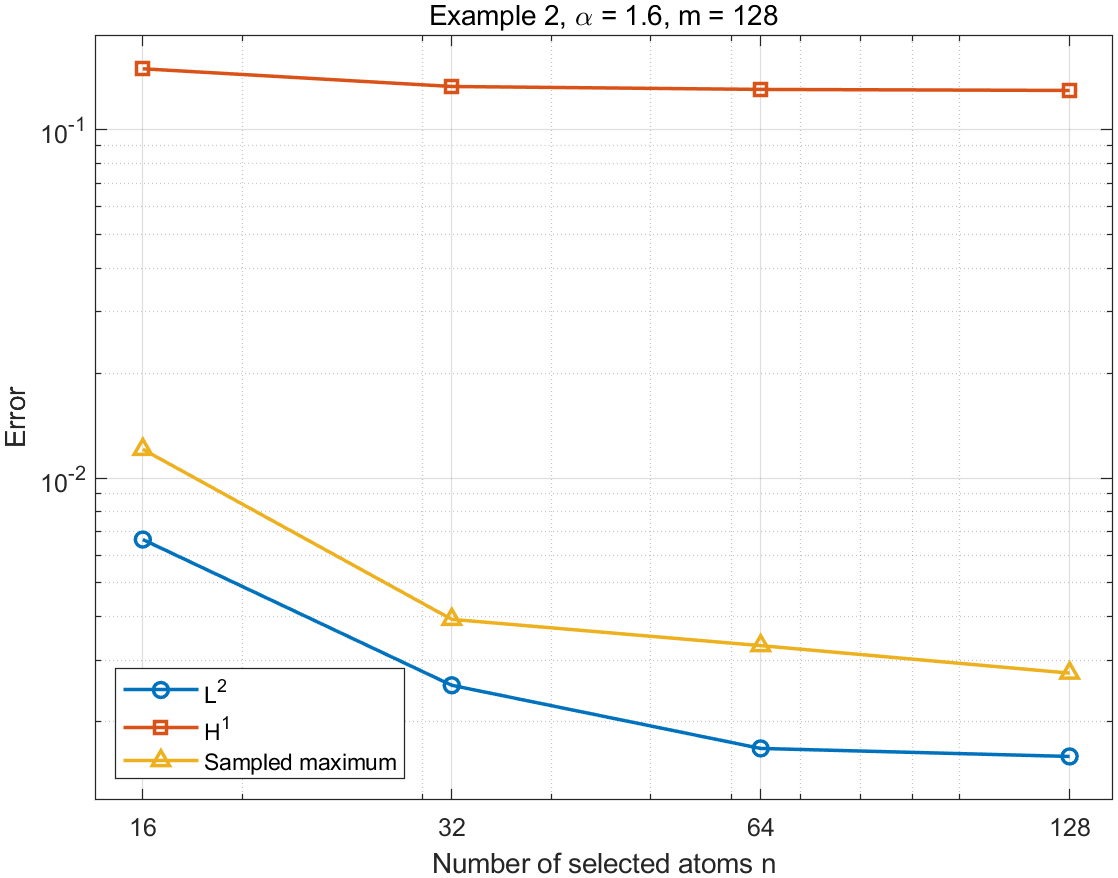}
\caption{Example 2: continuous error estimates against selected-atom count,
with $\alpha=1.2$ (left) and $\alpha=1.6$ (right). The maximum error is sampled.}
\label{fig:2d-ex2-convergence}
\end{figure}

\begin{figure}[htbp]
\centering
\includegraphics[width=0.84\textwidth]{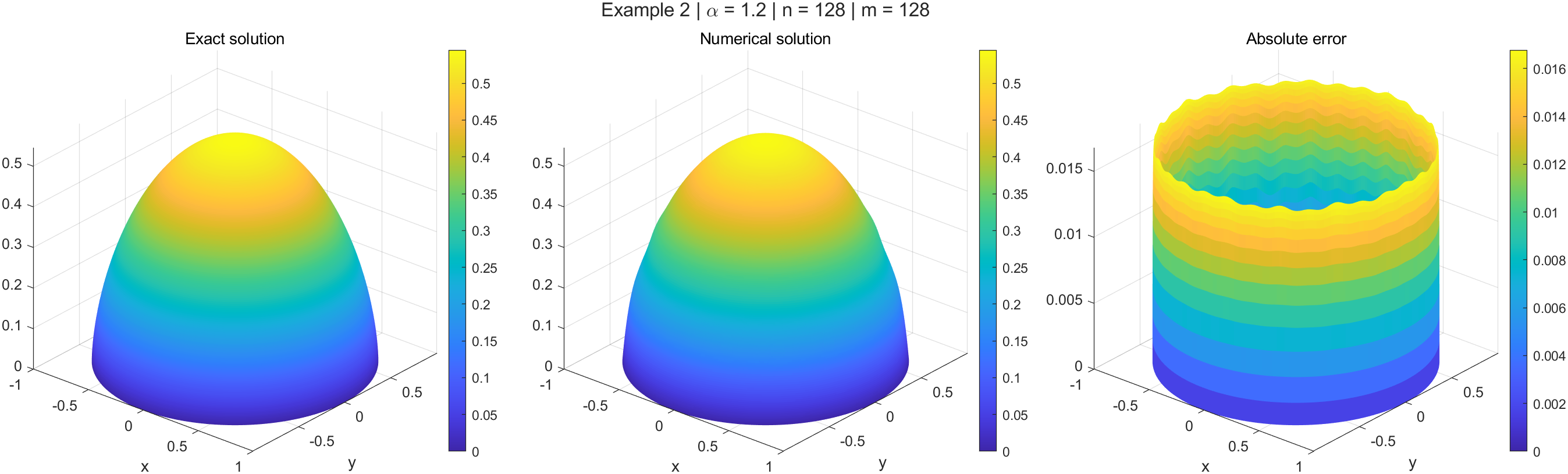}
\caption{Example 2, $\alpha=1.2$, $n=128$: exact solution (left),
network approximation (middle), and absolute error (right).
The first two panels share vertical and color scales.}
\label{fig:2d-ex2-a1p2-solutions}
\end{figure}

\begin{figure}[htbp]
\centering
\includegraphics[width=0.84\textwidth]{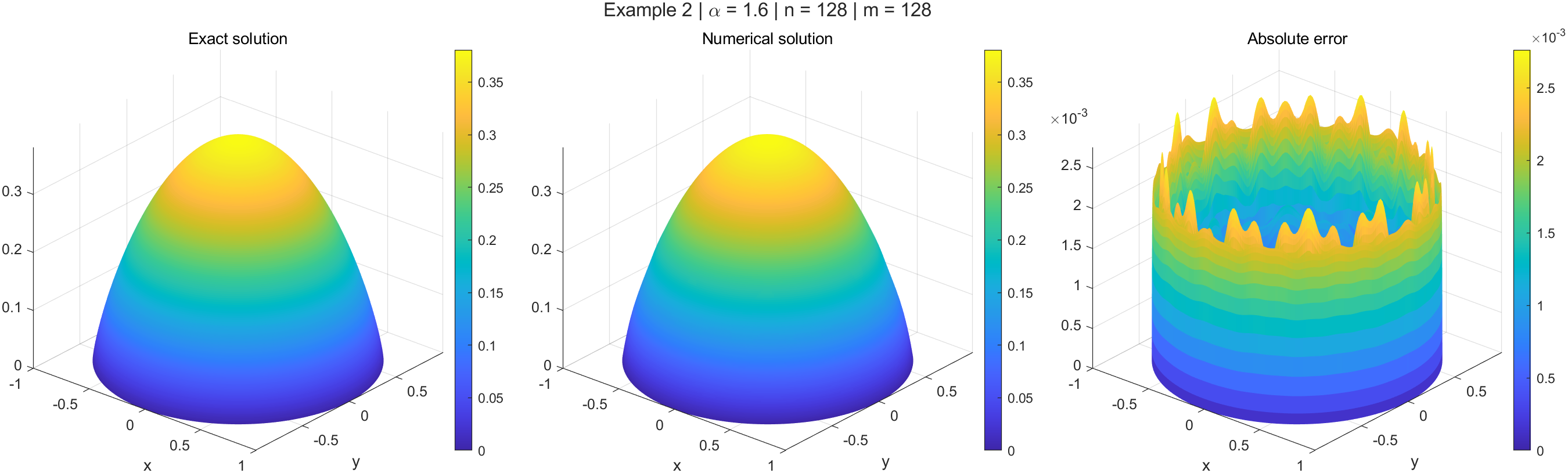}
\caption{Example 2, $\alpha=1.6$, $n=128$: exact solution (left),
network approximation (middle), and absolute error (right).
The first two panels share vertical and color scales.}
\label{fig:2d-ex2-a1p6-solutions}
\end{figure}

Figure~\ref{fig:2d-ex2-convergence} shows the pronounced slowing of the
$H^1$ improvement. Although the solution profiles in
Figures~\ref{fig:2d-ex2-a1p2-solutions} and
\ref{fig:2d-ex2-a1p6-solutions} are close in the interior, the error
surfaces expose a boundary layer. Each finite expansion in the present
dictionary has a bounded interior gradient, whereas the exact solution
has a singular gradient. This difference, together with the finite grid and dictionary resolution,
helps explain the slow $H^1$ improvement. The results motivate
boundary-adapted dictionary enrichment and adaptive resolution.

\FloatBarrier
As a check on the discretization, Table~\ref{tab:2d-grid-check} reports
the weighted nodal error $h\|R_hu-\boldsymbol u_h\|_2$ of the direct discrete
solution as the grid is refined, with $p_{\rm pad}=16$. The error decreases
in all four cases, more slowly for the boundary-singular solution.
Independent refinement of the integration and test grids changes the
final network errors by less than $0.005\%$ in $L^2$ and $H^1$ and
$0.37\%$ in the sampled maximum, supporting the reported error values.

\begin{table}[H]
\centering\small
\setlength{\tabcolsep}{5pt}
\caption{Grid errors of the direct discrete solution with $p_{\rm pad}=16$.
The four grids contain 193, 793, 3205 and 12849 interior nodes, respectively.}
\label{tab:2d-grid-check}
\begin{tabular}{c c r r r r}
\hline
Example & $\alpha$ & $m=16$ & $m=32$ & $m=64$ & $m=128$ \\
\hline
1 & 1.2 & $6.259\times10^{-3}$ & $1.563\times10^{-3}$ & $4.587\times10^{-4}$ & $1.249\times10^{-4}$ \\
1 & 1.6 & $8.324\times10^{-3}$ & $1.926\times10^{-3}$ & $5.151\times10^{-4}$ & $1.286\times10^{-4}$ \\
2 & 1.2 & $3.253\times10^{-2}$ & $1.985\times10^{-2}$ & $1.370\times10^{-2}$ & $7.122\times10^{-3}$ \\
2 & 1.6 & $3.006\times10^{-2}$ & $1.806\times10^{-2}$ & $1.090\times10^{-2}$ & $5.573\times10^{-3}$ \\
\hline
\end{tabular}
\end{table}
\FloatBarrier

\section{Conclusion}
We have extended energy-based orthogonal greedy approximation to
integral fractional Poisson equations by coupling OGA with discrete
fractional operators. Order-dependent GL matrices in one dimension and
a restricted fractional lattice operator on the disk provide positive
definite energy forms. Boundary-corrected and boundary-weighted neural
network dictionaries impose homogeneous exterior data, while normalized
residual selection and projection determine the approximation.

The discrete positivity and projection results explain how the OGA
framework carries over to these nonlocal problems. The error
decompositions distinguish the selected-space approximation from
operator discretization and, in two dimensions, finite-period effects.
The numerical experiments exhibit decreasing errors for the smooth
one-dimensional solution and the bounded-gradient disk benchmark.
The constant-forcing disk problem shows slower $H^1$ improvement,
highlighting the role of boundary regularity in dictionary design.

The method provides a computational extension of OGA to fractional
operators within a common energy framework. Further developments
include boundary-adapted dictionaries, joint refinement of the training
grid and periodic box, and FFT-based iterative solution strategies.
These directions can improve approximation near singular boundaries
and extend the approach to larger problems.

\section*{CRediT authorship contribution statement}
\textbf{Ruitong Shan:} Writing – review \& editing, Writing – original draft, Validation, Software, Methodology, Investigation. \textbf{Young Ju Lee:} Writing – review \& editing, Writing – original draft, Supervision, Methodology, Investigation, Funding acquisition, Conceptualization. \textbf{Jiwei Jia:} Writing – review \& editing, Writing – original draft, Supervision, Methodology, Investigation, Funding acquisition, Conceptualization.

\section*{Declaration of competing interest}
The authors declare that they have no known competing financial interests or personal relationships that could have appeared
to influence the work reported in this paper.

\section*{Data availability
}
Data will be made available on request.

\section*{Acknowledgments}
Authors would like to express their sincere appreciation to Professor Jinchao Xu for a number of his stimulating lectures on machine learning algorithms and guidance. The second author was funded in part by NSF-DMS 2208499 and the third author was funded in part by NSFC 22341302. 
\bibliographystyle{elsarticle-num-names} 
\bibliography{fref}
\end{document}